\documentclass[a4paper,12pt]{article}
\usepackage{times, url}
\usepackage{amsmath,amsthm} 
\usepackage[lite]{amsrefs}
\usepackage{chngcntr}
\usepackage[english]{babel}
\usepackage{enumitem}
\usepackage{hyperref} 
\usepackage{amsthm}
\usepackage{amssymb}
\usepackage[all,cmtip]{xy}
\usepackage{amsmath,amssymb,amsthm,amsfonts}

\newtheorem{theorem}{Theorem}[section]
\newtheorem{lemma}[theorem]{Lemma}

\theoremstyle{definition}
\newtheorem{definition}{Definition}[section]

\newtheorem{remark}{Remark}
\newtheorem*{case 1*}{Case 1}
\newtheorem*{case 2*}{Case 2}

\begin{document}

\setcounter{page}{1}

\begin{center}
{\LARGE \bf  On subset sums of $\mathbb{Z}_n^{\times}$ which are equally distributed modulo $n$.}
\vspace{8mm}

{\large \bf Gaitanas Konstantinos}
\vspace{3mm}

Department of  Applied Mathematical and Physical Sciences\\
 National Technical University of Athens \\ 
Heroon Polytechneiou Str., Zografou Campus, 15780 Athens, Greece \\ 
e-mail: \url{kostasgaitanas@gmail.com}
\vspace{2mm}

\end{center}
\vspace{10mm}

\noindent
{\bf Abstract:} In this paper we provide some results concerning the structure of a multiset $A$ with elements from $\mathbb{Z}_n^{\times}$, which has non-empty subset sums equally distributed modulo $n$. Here, $\mathbb{Z}_n^{\times}$ denotes the set which contains all the invertible elements of the ring $\mathbb{Z}_n$. In particular, we prove that if $n$ belongs to a certain subset of the natural numbers, then $A$ is a union of sets of the form $\{ a\cdot(\pm2^i)\}$. Additionally, we count the number of subsets of $\mathbb{Z}_q^{\times}$ with non-empty subset sums equally distributed modulo $q$, if $q$ is an odd prime power.
\vspace{10mm}

\section{Introduction} 
Let $A$ be a multiset of numbers. The sums representable in the form $\sum\limits_{a\in A}\epsilon_a\cdot a$ where $\epsilon_a\in \{0, 1\}$ are called the \emph{subset sums} of $A$ and the multiset containing \emph{all} the subset sums is often denoted by $\sum(A)$ that is, 
\begin{gather*}\sum(A)=\{\sum_{x\in B}x: B\subseteq A\}\end{gather*} Let $\sum(A)^{(s)}$ denote the multiset which contains the $\binom{|A|}{s}$ subset sums of $A$ which have size $s$.
Through the last years there has been a growing interest in the following problem: 
\begin{center}\emph{Is it possible to recover $A$, given $\sum(A)^{(s)}$?} \end{center}
This problem is credited to Leo Moser, who proposed to recover $A$, given it's $\binom{|A|}{2}$ $2$-subset sums, if $|A|=5$. The general case of this problem was first studied by J. L. Selfridge and E. G. Straus in \cite{SelStra}. The most recent progress on this problem can be found in \cite{AnFe}, where Andrea Ciprietti and Federico Glaudo studied this problem from another point of view; they proved that for an abelian group $A$, the reconstruction problem is solvable if and only if every order of a torsion element of the abelian group satisfies a certain number-theoretical property. For more related material that might satisfy the reader's curiosity we refer to the survey of Dmitri V. Fomin \cite{Fom}. In this paper, our main objective is to explicitly determine $A$, on condition that equidistribution modulo $n$ occurs among it's non-emtpy subset sums. 
\subsection{A preliminary theorem}
We begin by giving a proof that the non-empty subset sums of $\mathbb{Z}_n^{\times}$ are equally distributed modulo $n$, if $n$ is odd (if $n$ is even, this is not possible, as we will see). It seems a near certainty that a proof of this result has already appeared somewhere in the long history of mathematics, but I have not yet been able to find a reference. Thereby, we provide a simple proof which does not seem to appear in the literature and hope also that will be accessible to a broad spectrum of readers. Furthermore, the primary motivation for the arguments we present in this paper arises from the construction presented in the proof. In the following, $\varphi(n)$ denotes Euler's totient function which counts the number of positive integers not greater than $n$, which are relatively prime to $n$. We also define the order of $2$ modulo $n$, as the least exponent $r$ for which the congruence $2^r\equiv 1\pmod{n}$ holds. Note that from Euler's theorem we know that $r\mid \varphi(n)$.
\begin{theorem}\label{th} Let $n\in \mathbb{N}$ be odd. The non-empty subset sums of $\mathbb{Z}_n^{\times}$ are equally distributed modulo $n$. In particular, every residue class appears exactly $\frac{2^{\varphi(n)}-1}{n}$ times.
\begin{proof}We slightly modify the construction presented by the author in \cite{Gaitanas}. Let $r$ be the order of $2$ modulo $n$ and $k=\frac{\varphi(n)}{r}$. We consider the following subsets of $\mathbb{Z}_n^{\times}$: For $i=1, \ldots, k$, we let \begin{gather*}B_i=\{b_i\cdot 2^{m-1} \pmod{n}: 1\leq m\leq r\}\end{gather*} with the restriction that $(b_i, n)=1$ and $b_i\not \equiv b_{\kappa}\cdot2^{m-1}\pmod{n}$ for $\kappa<i$.\\
Observe that the subsets $B_i$ do not have any element in common. To see this, assume for the sake of contradiction that $B_i\cap B_j\neq \emptyset$ for some $i, j$ not greater than $k$. Then there is an $m'\leq r$, such that $b_i\cdot 2^{m-1}\equiv b_j\cdot 2^{m'-1}\pmod{n}$. From this we obtain that $b_i \equiv b_j\cdot 2^{m'-m}\pmod{n}$ which is absurd. Thus, the union of $B_i$ contains $k\cdot r=\varphi(n)$ distinct elements, which is equivalent to \begin{gather*}B_1\cup B_2 \cup \cdots\cup B_k=\mathbb{Z}_n^{\times}.\end{gather*}
In addition, it is evident that the non-empty subset sums of $B_i$ produce every element of the set $X_i=\{b_i\cdot x:  1\leq x\leq 2^r-1\}$, since every binary representation of $x$ occurs exactly once. Since $n\mid 2^r-1$ and $(b_i, n)=1$, the residue classes $\pmod{n}$ are equally distributed (they occur exactly $\frac{2^r-1}{n}$ times) in $X_i$ and so do in $\sum(B_i)$. It is relatively easy to see that their union also produces equally distributed subset sums and thus every residue class is represented  $\frac{2^{\varphi(n)}-1}{n}$ times. This proves the theorem.\footnote{This method was used also by Euler himself to prove that $a^{\varphi(n)}\equiv 1\pmod{n}$.}
\end{proof}
\end{theorem}
Inspired by the method of the previous proof, we ask the following (converse) question: If the non-empty subset sums of $A\subseteq\mathbb{Z}_n^{\times}$ are equally distributed modulo $n$, then is $A$ necessarily the union of sets consisting of elements which are in ``geometric progression" with ratio $2$? We show that for a certain subset of $\mathbb{N}$ this is almost true. In particular, the main result of this paper is that the ratio must be $\pm2\pmod{n}$. 

  \section{Main results}
In this section we state our main results and provide some key proof techniques and insights. Before we begin, we would like to mention that in the following, Theorems \ref{theore} and \ref{pr} are very similar to proposition 4.1 in \cite{AnFe}. We begin with the following.
\begin{theorem}\label{theore}
Let $A$ be a multiset with elements from $\mathbb{N}$, $|A|=k$ and $r$ be the order of $2\pmod{n}, n\geq3$. Suppose that the non-empty subset sums of $A$ are equally distributed $\pmod{n}$. Then the following hold true:
\counterwithin{enumi}{theorem}
\begin{enumerate}
\item\label{part1}$\prod_{i=1}^{k} (1+x^{a_i})=1+R(x)\cdot \frac{x^n-1}{x-1}$, where $R(x)\in \mathbb{Z}[x]$.
\item\label{part2}$2^{k}\equiv 1\pmod{n}$ (which implies that $n$ is odd and $r\mid k$).
\item\label{part3}$\sum_{i=1}^{k}a_i\equiv 0\pmod{n}$.
\end{enumerate}
\begin{proof} It is relatively easy to prove \ref{part1} and \ref{part2}. First, we expand the product $\prod_{i=1}^{k} (1+x^{a_i})$ as a polynomial and take every summand modulo $\frac{x^n-1}{x-1}$. If we denote by $A_m$ the number of non-empty subset sums which are congruent to $m\pmod{n}$, we obtain\footnote{Here we adopt the convention that the empty sum is $0\pmod{n}$} \begin{gather*}\prod_{i=1}^{k} (1+x^{a_i})=1+A_0x^0+A_1x^1+\ldots+ A_{n-1}x^{n-1}+Q(x)\cdot \frac{x^n-1}{x-1},\end{gather*} where $Q(x)\in \mathbb{Z}[x]$. But the subset sums are equidistributed, so $A_m=\frac{2^{k}-1}{n}$ for every $m$. Thus, \begin{gather*}\prod_{i=1}^{k} (1+x^{a_i})=1+\frac{2^{k}-1}{n}(x^0+\ldots+x^{n-1})+Q(x)\cdot \frac{x^n-1}{x-1},\end{gather*} which proves \ref{part1} with $R(x)=\frac{2^{k}-1}{n}+Q(x)$.\\
Now we take the limit as $x$ approaches $1$ on both sides of \ref{part1} and obtain $2^{k}=1+R(1)\cdot n$. Since $R(1)\in \mathbb{Z}$, we have $2^{k}\equiv 1\pmod{n}$ which proves \ref{part2}.\\
To see that $\sum_{i=1}^{k} a_i\equiv 0\pmod{n}$, we differentiate both sides of \ref{part1} and get \begin{gather*}\prod_{i=1}^{k} (1+x^{a_i})\cdot \sum_{i=1}^{k}\frac{a_ix^{a_i-1}}{1+x^{a_i}}=R'(x)\cdot \frac{x^n-1}{x-1}+R(x)\cdot \frac{(n-1)x^n-nx^{n-1}+1}{(x-1)^2}.\end{gather*} Again, we take the limit of both sides as $x$ approaches $1$. We obtain that \begin{gather*}2^{k-1}\cdot \sum_{i=1}^{k} a_i=R'(1)\cdot n+R(1)\frac{n(n-1)}{2}.\end{gather*} Since $n$ is odd, the right-hand side is divisible by $n$ and so is the left-hand side. This yields $2^{k-1}\cdot \sum_{i=1}^{k}a_i\equiv0\pmod{n}$ and thus $\sum_{i=1}^{k} a_i\equiv0\pmod{n}$. With this, we proved \ref{part3}.\\
\end{proof}
\end{theorem}
Before we state our main theorem, we provide an auxiliary result. It appears in \cite{Cha&Dhi} and is presented as two separate theorems. For the convenience of the reader, we present them together and slightly reshape them under the restriction that we consider only odd numbers. We say that $m$ is a semi-primitive root $\pmod{q}$ if the order of $n\pmod{q}$ is $\frac{\varphi(q)}{2}$ and with $\bar{\mathbb{Q}}$ we denote the field of algebraic numbers.
\begin{theorem}(Chatterjee, T., Dhillon, S. \cite{Cha&Dhi})\label{ChaDhi}
\counterwithin{enumi}{theorem}
\begin{enumerate}
 \item\label{part 4} Let $q\geq 3$ be an odd integer. The elements of \begin{gather*} L_1=\{\log(2\sin\frac{a\pi}{q}): 1\leq a\leq \frac{q-1}{2},\quad (a, q)=1\}\quad \text{and}\quad \pi\end{gather*} are linearly independent over $\bar{\mathbb{Q}}$ if and only if $q$ is a prime power.
 \item\label{part 5} Let $q\geq 3$ be an odd integer with at least two prime divisors. The elements of \begin{gather*}L_2=\{\log(2\sin\frac{a\pi}{q}): 1<a\leq \frac{q-1}{2},\quad (a, q)=1\} \quad \text{and}\quad \pi\end{gather*} are linearly independent over $\bar{\mathbb{Q}}$ if and only if $q$ satisfies one of the following conditions: (here $a_1, a_2, a_3\geq 1$ and $p_1, p_2, p_3$ are odd primes)
\begin{enumerate}
\item[\textbf{(I)}] $q=p_1^{a_1}p_2^{a_2}$; and
\item[(Ia)] when $p_1\equiv p_2\equiv 3\pmod{4}$: $p_1$ is a semi-primitive root $\pmod{p_2^{a_2}}$ and $p_2$ is a semi-primitive root $\pmod{p_1^{a_1}}$ or vice versa. 
\item[(Ib)]Otherwise: $p_1$ and $p_2$ are primitive roots $\pmod{p_2^{a_2}}$ and $\pmod{p_1^{a_1}}$ respectively.
\item[\textbf{(II)}] $q=p_1^{a_1}p_2^{a_2}p_3^{a_3}$; $p_1\equiv p_2\equiv p_3\equiv 3\pmod{4}$: $\frac{p^i-1}{2}, (1\leq i\leq 3)$ are co-prime to each other; and
\item[(IIa)] $p_1, p_2, p_3$ are primitive roots $\pmod{p_2^{a_2}}$, $\pmod {p_3^{a_3}}$, $\pmod{p_1^{a_1}}$, respectively and semi-primitive roots $\pmod{p_3^{a_3}}$, $\pmod {p_1^{a_1}}$, $\pmod{p_2^{a_2}}$, respectively.
\end{enumerate}
\end{enumerate}
\end{theorem}
In the original paper the number $\ln 2$ appears among the elements of $L_1, L_2$ and they also require $\sin(\frac{a\pi}{q})\neq \frac{1}{2^{\alpha}}, \alpha \in \mathbb{Q}$. This version of the theorem appears in order to avoid the trivial case where $\frac{a\pi}{q}=\frac{\pi}{6}$ which implies $\log(2\sin\frac{a\pi}{q})=0$. Since we consider only odd numbers $q$, we disregard this case. We will also make use of the following definition:
\begin{definition}We say that a set $B\subseteq\mathbb{Z}_n^{\times}$ is a geometric progression with ratio $l$ with leader $b$, if every element in $B$ is of the form $b\cdot l^{j-1}$ (modulo $n$).
\end{definition}
 Our main result is the following theorem.
\begin{theorem}\label{mtheor}
Let $A$ be a multiset with elements from $\mathbb{Z}_q^{\times}$, $|A|=k$ and $r$ be the order of $2\pmod{q}$. Suppose that the non-empty subset sums of $A$ are equally distributed. If
\begin{enumerate}
\item[(1)] $q$ is an odd prime power or
\item[(2)] $q$ satisfies the conditions of Theorem \ref{part 5} and $1,\frac{q-1}{2}, \frac{q+1}{2}, q-1 \not \in A$ 
\end{enumerate}
Then, there are disjoint sets $B_1, \ldots ,B_t\subseteq A$, such that $A=\cup_{i=1}^{t}B_i$ and each $B_i$ is a geometric progression with leader $b_i$ and ratio $\pm2$. Additionally, $|B_i|=r$, for every $i\leq t$.
\begin{proof}
We begin by proving the first part of the theorem. Let $\omega=e^{2\pi i/q}$ be a $q$-th root of unity where $q$ is a power of an odd prime. From Theorem \ref{part1} we know that  $\prod_{i=1}^{k} (1+x^{a_i})=1+R(x)\cdot \frac{x^q-1}{x-1}$. For $x=\omega$, we derive that $\prod_{i=1}^{k}(1+\omega^{a_i})=1$. But, \begin{gather*}1+\omega^{a_i}=2\cos\frac{\pi a_i}{q}\cdot e^{i\cdot \pi a_i/q}.\end{gather*} It follows that \begin{gather*}2^k\prod_ {i=1}^{k}\cos \frac{\pi a_i}{q}\cdot e^{\pi\cdot i\sum a_i/q}=1.\end{gather*} We multiply both sides by $2^k\cdot \prod_{i=1}^{k}\sin \frac{\pi a_i}{q}$ and use the identity $\sin 2x=2\cos x\cdot \sin x$. Also note that from Theorem \ref{part3}, $q\mid \sum_{i=1}^{k}a_i$, which implies that $ e^{\frac{{\sum_{i=1}^{k}a_i\cdot \pi i}}{q}}=\pm1$. We obtain \begin{equation}\label{eq:1}\prod_{i=1}^{k}\left (2\sin \frac{2\pi a_i}{q}\right )=\pm \prod_{i=1}^{k}\left(2\sin \frac{\pi a_i}{q}\right ).\end{equation}
Now define \begin{gather*}\langle x\rangle=\min\{x\pmod{q},\quad q-x\pmod{q}\}.\end{gather*} This is the residue of $x$ modulo $q$ which has absolute value less than $q/2$. It will be useful to note that $\langle x\rangle=\langle y\rangle$ implies that $x\equiv \pm y\pmod{q}$. Thus, by observing that $\frac{\pi \langle x\rangle}{q}<\frac{\pi}{2}$ and $\sin \frac{\pi x}{q}=\pm \sin \frac{\pi \langle x\rangle}{q}$, we can rewrite equation \ref{eq:1} as \begin{equation}\label{eq:2}\prod_{i=1}^{k}\left( 2\sin \frac{\langle 2a_i\rangle\pi}{q} \right )=\prod_{i=1}^{k}\left (2\sin \frac{\langle a_i\rangle\pi}{q}\right ).\end{equation} We apply logarithms to both sides (we may take the principal branch of the logarithm) and get:\begin{equation} \label{eq:3}\sum_{i=1}^{k}\log\left(2\sin\frac{\langle 2a_i\rangle\pi}{q}\right)=\sum_{i=1}^{k}\log\left(2\sin\frac{\langle a_i\rangle\pi}{q}\right)+M\pi\cdot i, M\in \mathbb{Z}.\end{equation}
\footnote{Here $M\pi\cdot i$ denotes the imaginary unit multiplied by an integer multiple of $\pi$.} Notice that for every $i\leq k$, both $\langle a_i\rangle$ and $\langle 2a_i\rangle$ are co-prime with $q$ and not greater than $\frac{q-1}{2}$. Thus, we may apply Theorem \ref{ChaDhi}. and immediately get that $M=0$. Using the notation of Theorem \ref{part 4}, equation \ref{eq:3} can be written in the form \begin{gather*}\sum_{x_i\in L_1} x_i=\sum_{y_i\in L_1} y_i.\end{gather*} But $L_1$ is a linearly independent set. Hence every summand on the left-hand side of the equation appears exactly the same amount of times at the right-hand side. This means that we are able to attach to every $i$ a unique $s$, such that $\log\left(2\sin\frac{\langle 2a_i\rangle\pi}{q}\right)=\log\left(2\sin\frac{\langle a_s\rangle\pi}{q}\right)$. The quantities contained in the logarithms are positive and it is evident that $0<\frac{\langle 2a_i\rangle\pi}{q},\frac{\langle a_s\rangle\pi}{q}<\frac{\pi}{2}$. We conclude that $\langle 2a_i\rangle=\langle a_s\rangle$, which implies \begin{gather*}a_s\equiv\pm 2a_i\pmod{q}.\end{gather*} Now there is a simple path to finish the proof: We relabel the elements of $A$. We choose $a_1$ to be an arbitrary element from $A$ and proceed by choosing $a_j$ to be the (unique) element such that $a_j\equiv\pm2a_{j-1}$ for $j=2,...,r$. This is equivalent to \begin{gather*}a_j\equiv a_1\cdot(\pm2^{j-1})\pmod{q} \quad , \quad 1\leq j\leq r.\end{gather*} From this construction it is also evident that $a_1\equiv\pm 2a_r$ since $2^r\equiv 1\pmod{q}$. We call this set $B_1=\{a_1,\ldots ,a_r\}=\{b_1\cdot (\pm2^{j-1})\}$ with leader $b_1=a_1$. We can repeat the same construction (to the multiset $A\setminus B_1$) and eventually we will exhaust all the elements of $A$. That settles the theorem if $q$ is a power of an odd prime.\\ 
We now focus on the second case. In order to not be too repetitive we only sketch briefly what modifications are to be made. We may repeat exactly the same arguments as long as $\log\left(2\sin\frac{1\cdot \pi}{q}\right)$ does not emerge in any of the sides of equation \ref{eq:2}. It is a routine matter to see that $\langle a_i\rangle=1$ only if $a_i$ equals to $1$ or $q-1$. Also, $\langle 2a_i\rangle=1$ only if $a_i$ equals to $\frac{q-1}{2}$ or $\frac{q+1}{2}$. But this is impossible, since from our assumption $\{1, \frac{q-1}{2}, \frac{q+1}{2}, q-1\}\cap A=\emptyset$. This completes the proof.
\end{proof}
\end{theorem}
\begin{remark} For a moment, we take a brief look at the well-known trigonometric identity \footnote{This is often called ``Morrie's Law" \cite {Morrie}.}\begin{equation}\label{eq:4}2^k\cos(a)\cdot \cos(2a)\cdots \cos(2^{k-1}a)=\frac{\sin(2^ka)}{\sin a}.\end{equation}
One may ask if the angles must be in a geometric progression with ratio $2$ in order to have equality on both sides. In particular, we may consider the case where there are $k$ angles which are rational multiples of $\pi$ with denominator $q$. Suppose that \begin{gather*} 2^k\cos \left(\frac{\pi a_1}{q}\right)\cdot \cos \left(\frac{\pi a_2}{q}\right)\cdots \cos \left(\frac{\pi a_k}{q}\right)=\frac{\sin \left(2^k\frac{\pi a_1}{q}\right)}{\sin \left(\frac{\pi a_1}{q}\right)}.\end{gather*} Since $2^k\equiv 1\pmod{q}$, this transposes to \begin{gather*} \prod_{i=1}^{k}\sin \left(\frac{2\pi a_i}{q}\right)=\pm \prod_{i=1}^{k}\sin \left(\frac{\pi a_i}{q}\right)\Rightarrow \prod_{i=1}^{k}\sin \left(\frac{\langle 2a_i\rangle\pi}{q}\right)=\prod_{i=1}^{k}\sin \left(\frac{\langle a_i\rangle\pi}{q}\right).\end{gather*} As we have seen, this is possible only if all the $a_i$'s belong to union of sets where each one is a geometric progression with ratio $\pm 2$. Thus, as a corollary of our main theorem we see (informally) that identity \ref{eq:4} holds true only for a ``unique'' choice of angles.
\end{remark} 
In the sequel we shall prove the following theorem, which can be regarded as the converse of Theorem \ref{mtheor}. Although we will eventually focus in the case where $q$ is a prime power, it is worth mentioning that this result holds true for every odd number $n$.
\begin{theorem}\label{pr}Let $n$ be odd and $2^r\equiv 1\pmod{n}$. Suppose that $A$ is a multiset with elements from $\mathbb{Z}_n$ and $|A|=k$. If $A=\cup_{i=1}^{t}B_i$, where $B_i=\{b_i\cdot (\pm2^{j-1}): 1\leq j \leq r\}$ and $\sum_{a\in A}a\equiv 0\pmod{n}$, the non-empty subset sums of $A$ are equally distributed modulo $n$.
\begin{proof}
We denote by $S^+$ the sum of all the elements from the sets $B_1, \ldots, B_t$ which are of the form $b_i(+2^{j-1})$ and $S^-$ the sum of all these which are of the form by $b_i(-2^{j-1})$. Additionally, let $A_m$ be the number of non-empty subset sums which are congruent to $m\pmod{n}$.
For a moment, we concentrate on the product \begin{equation}\label{eq:5}\mathcal{P}_x=\prod_{a\in A}(1+\omega^{x\cdot a})=1+\sum_{m=0}^{n-1}A_m\omega^{x\cdot m}.\end{equation} Here $\omega=e^{2\pi i/n}$ is a primitive $n$-th root of unity. It is straightforward to see that for $x=0$, we have $\mathcal{P}_0=2^{k}$. For $x\geq 1$, since $B_i=\{b_i\cdot (\pm2^{j-1})\}$ for $1\leq i\leq t$, we may write \begin{gather*}\mathcal{P}_x=\prod_{i=1}^{t}\prod_{j=1}^{r}\left(1+\omega^{x\cdot b_i\cdot(\pm2^{j-1})}\right)=\omega^{x\cdot S^-}\prod_{i=1}^{t}\prod_{j=1}^{r}(1+\omega^{x\cdot b_i\cdot2^{j-1}}).\end{gather*}
But \begin{gather*}\prod_{j=1}^{r}(1+\omega^{x\cdot b_i\cdot2^{j-1}})=\frac{1-\omega^{x\cdot b_i\cdot2^r}}{1-\omega^{x\cdot b_i}}=1,\end{gather*} since $2^r\equiv1\pmod{n}$ and $\omega^{x\cdot b_i}\neq 1$. Thus, \begin{gather*}\mathcal{P}_x=\omega^{x\cdot S^-}, \quad x\geq 1.\end{gather*}
From equation \ref{eq:5}, it is possible to obtain the inverse Fourier transform of $\{A_m\}_{m=0}^{n-1}$. From the inverse transform formula, we obtain that for $m\geq 1$ \begin{gather*}A_m=\frac{1}{n}\sum_{x=0}^{n-1}\mathcal{P}_x\cdot \omega^{-x\cdot m}=\frac{1}{n}(2^{k}+\sum_{x=1}^{n-1}\omega^{x\cdot(S^--m)}).\end{gather*}
If $m=0$, $A_0$ is equal to the same value minus $1$ since we adopt the convention that the empty sum is $0\pmod{n}$. Hence \begin{gather*}A_0=-1+\frac{1}{n}(2^{k}+\sum_{x=1}^{n-1}\omega^{x\cdot(S^--m)}).\end{gather*}We observe that $\sum_{x=1}^{n-1}\omega^{x\cdot(S^--m)})$ is a sum of $n-1$ terms who belong to a geometric progression. It is well-known that \begin{gather*}\sum_{x=1}^{n}\omega^{x\cdot c}=\begin{cases} 
      n, & \quad\text{if}\quad  n\mid c \\
      0, &  \quad\text{if}\quad n\nmid c\\
   \end{cases}\end{gather*}  It follows that

\begin{gather*}\sum_{x=1}^{n-1}\omega^{x\cdot(S^--m)})=\begin{cases} 
      n-1, & \quad\text{if}\quad  n\mid S^--m \\
      -1, &  \quad\text{if}\quad n\nmid S^--m\\
   \end{cases}\end{gather*}
Putting this all together we see that for $m\neq 0$,

\begin{gather*}A_m=\begin{cases} 
      \frac{1}{n}(2^{k}-1) & \quad\text{if}\quad  S^-\not\equiv m\pmod{n} \\
       \frac{1}{n}(2^{k}-1)+1 &  \quad\text{if}\quad S^-\equiv m\pmod{n}\\
   \end{cases}\end{gather*}
Also, for $m=0$,
\begin{gather*}A_0=\begin{cases} 
      \frac{1}{n}(2^{k}-1)-1 & \quad\text{if}\quad  S^-\not\equiv 0\pmod{n} \\
       \frac{1}{n}(2^{k}-1) &  \quad\text{if}\quad S^-\equiv 0\pmod{n}\\
   \end{cases}\end{gather*} 
Since $\sum_{a\in A}a\equiv0\pmod{n}$, we have $S^++S^-\equiv 0\pmod{n}$. Also, it is evident that $S^+-S^-\equiv 0\pmod{n}$ since this is equivalent to \begin{gather*}b_1(1+2^1+\ldots+2^{r-1})+\ldots+b_t(1+2^1+\ldots+2^{r-1})\equiv0\pmod{n}.\end{gather*} This is obviously true because every summand  is congruent to $0$ modulo $n$. From this, we obtain that $S^+\equiv S^-\equiv0\pmod{n}$.
Since $S^-\equiv0\pmod{n}$, we obtain $A_0=A_1=\ldots= A_{n-1}= \frac{1}{n}(2^{k}-1)$.
\end{proof}
\end{theorem}
\begin{remark}
We pause now to record an observation that follows readily from the computations made above. If we omit the restriction $\sum_{a\in A}a\equiv 0\pmod{n}$, among the non-empty subset sums equidistribution modulo $n$ occurs except possibly for only two residue classes. In particular, these two residue classes differ from the others only by $1$. This is somewhat surprising, since we may multiply randomly by $+2^{j-1}$ or $-2^{j-1}$ every leader in $B_1, \ldots, B_t$ and still preserve rigid structure among the subset sums of $A$.
\end{remark} 
\section{Counting sets}
 Let $q$ be a power of an odd prime and $r$ be the order of $2$ modulo $q$. In this section, we will count the number $\mathcal{E}(q)$ of sets with elements from $\mathbb{Z}_q^{\times}$, having non-empty subset sums equidistributed modulo $q$. From Theorem \ref{part3}, Theorem\ref{mtheor} and Theorem \ref{pr}, this is possible if and only if $A$ is a union of sets of the form \begin{gather*}B_i=\{b_i\cdot (\pm2^{j-1}): 1\leq j \leq r\},\end{gather*} such that $\sum_{a\in A}a\equiv0\pmod{q}$. This forces us to put some restraint in the way we define the leaders $b_i$ and the choice of plus or minus, in order to have distinct elements which in addition sum to $0$ modulo $q$. For the reader’s convenience, we will shortly explain the general strategy.
\subsection{The number of possible leaders.} For the rest of the paper, we will denote the number of possible leaders with $t$. We begin by making a construction similar to the one presented in Theorem \ref{th}. We choose the first leader to be $b_1=1$ and define the $i$-th leader $b_i$ as the least element in $\mathbb{Z}_q^{\times}$ which has the property that \begin{gather*}b_i\not \equiv b_{\kappa}(\pm2^{j-1})\pmod{q}:\quad {\kappa}<i, \quad j=1,\ldots, r.\end{gather*} In order to count the number of possible leaders, we have to treat separately the case where $r$ is odd or even. We recall that the number of elements of $\mathbb{Z}_q^{\times}$ is $\varphi(q)$.
\begin{case 1*}\emph{The order of $2$ modulo $q$ is odd.}\\ 
If $r$ is odd, for every possible choice of plus or minus the elements of $\{ b_i(\pm 2^{j-1}): j=1,\ldots, r\}$ are incogruent modulo $q$. To see this, suppose that $b_i(+2^{j-1})\equiv b_i(+2^{w-1})\pmod{q}$ or $b_i(+2^{j-1})\equiv b_i(-2^{w-1})\pmod{q}$ for some $j, w\leq r$. Since $(b_i, q)=1$, the first implies that $2^{j-w}\equiv 1\pmod{q}$ and the second that $2^{j-w}\equiv-1\pmod{q}$. We conclude that $r\mid j-w$ or $r\mid2(j-w)$. Both are impossible, since $j-w$ is less than $r$ and $r$ is odd. Thus, for all possible choices of plus and minus, all numbers of the form $b_i(\pm 2^{j-1})$ span $2r$ elements of $\mathbb{Z}_q^{\times}$. From this and the definition of $b_i$, we obtain that there are $t=\frac{\varphi(q)}{2r}$ possible leaders.
\end{case 1*}
\begin{case 2*}\emph{The order of $2$ modulo $q$ is even.}\\ 
If $r$ is even, from $2^r\equiv 1\pmod{q}$ we obtain $(2^{r/2}-1)(2^{r/2}+1)\equiv0\pmod{q}$. But $2^{r/2}-1$ and $2^{r/2}+1$ are coprime, since they are both consecutive odd numbers. It follows that $q\mid 2^{r/2}-1$, or $q\mid 2^{r/2}+1$. Since $r$ is the order of $2$ modulo $q$ we must have $2^{r/2}\equiv-1\pmod{q}$. Hence, $b_i\cdot 2^{j-1}$ and $b_i(-2^{j-1+r/2})$ are congruent, which implies they cannot occur simultaneously in $B_i$. In particular, both $b_i\cdot 2^{j-1}$ and $b_i\cdot2^{j-1+r/2}$ must have the same sign in $B_i$ for every $j=1, \ldots, r/2$. Thus, for all possible choices of plus and minus, all numbers of the form $b_i(\pm 2^{j-1})$ span $r$ elements of $\mathbb{Z}_q^{\times}$. From this and the definition of $b_i$, we obtain that there are $t=\frac{\varphi(q)}{r}$ possible leaders.\end{case 2*}
\begin{remark}From the previous analysis, the careful reader may notice that every set\\
 $B=\{b\cdot (\pm2^{j-1}): 1\leq j \leq r\}$ is equal to some set $B_i=\{b_i\cdot (\pm2^{j-1}): 1\leq j \leq r\}$. This is trivially true if $b=b_i$, for some leader $b_i$. If $b$ is not equal to any of the leaders, then from definition we must have $b\equiv b_i(\pm 2^{x-1})\pmod{q}$ for some leader $b_i$ and for some $x\leq r$. Thus, the set $\{b\cdot (\pm2^{j-1}): 1\leq j \leq r\}$ is equal to the set $\{b_i\cdot (\pm2^{x+j-2}): 1\leq j \leq r\}$. Since $\{2^n\pmod{q}: n\in \mathbb{Z}\}$ is a cyclic group with order $r$,  the last set is eventually equal to the set $B_i=\{b_i\cdot (\pm2^{j-1}): 1\leq j \leq r\}$.\end{remark}
\subsection{The choice of ``$+$'' and ``$-$''.} What remains, is to count the number of ways we can choose plus or minus inside the sets $B_i$, in order to ensure that $\sum_{a\in A}a\equiv0\pmod{q}$. Recall that in Theorem \ref{pr} we denoted by $S^+$ the sum of all the elements from the sets $B_i$ which are of the form $b_i(+2^{j-1})$, by $S^-$ the sum of all these which are of the form $b_i(-2^{j-1})$ and proved that $S^{+}\equiv S^{-}\equiv 0\pmod{q}$. Thus, it will be enough to restrict ourselves to the choice of ``$+$''. We denote by $S_i^{+}$ the sum of all the elements of $B_i$ which are multiplied by $+2^{j-1}$. Suppose $A$ is a union of $\tau$ sets, that is,  $A=\cup_{i=1}^{\tau}B_i$. Thus, $\sum_{a\in A}a\equiv0\pmod{q}$ is equivalent to \begin{gather*}\sum_{i=1}^{\tau}S_i^{+}\equiv0\pmod{q}. \end{gather*}We observe that a choice of ``$+$'' inside $B_i$, corresponds to a sum of powers of $2$ all of which are not greater than $2^{r-1}$. Every natural number has a unique binary representation, hence a choice of ``$+$'' corresponds to a unique value of $S_i^{+}$. Thus, $S_i^{+}$ can take $2^r$ possible values. (Note that we want to take account the case where all the signs in $B_i$ are minus, which corresponds to $S_i^{+}=0$). Hence, $\mathcal{E}(q)$ is equal to the number of solutions of the following congruence:\begin{equation}\label{eq:6}\sum_{i=1}^{\tau} b_i\cdot x_i\equiv0\pmod{q}:\quad 0\leq x_i\leq2^r-1. \end{equation}
Thereby, we introduce an elementary lemma which as we will see, simplifies the proof of Theorem \ref{count}. In the following, $\mathbf{b}\cdot \mathbf{x}$ stands for the dot product of two vectors $\mathbf{b}=\left[b_1, \ldots, b_s\right]$ and $\mathbf{x}=[x_1, \ldots, x_s]$, defined by $\mathbf{b}\cdot \mathbf{x}=b_1x_1+\ldots+b_sx_s$. We will also use the term ``$N$- vector'' for a vector which has no repeated entries. 
\begin{lemma}\label{lemma}
   Suppose $B=\{b_1, \ldots, b_t\}\subseteq \mathbb{Z}_q^{\times}$ and $X=\{0, 1, 2, \ldots , kq\}$, $k\in \mathbb{N}^*$. Let $\mathcal{H}(n, \mathbf{b}, \mathbf{x})$ be the number of solutions of the congruence \begin{gather*}\mathbf{b}\cdot \mathbf{x}\equiv0\pmod{q},\end{gather*} where $\mathbf{b}$ is an $N$-vector with $n$ entries from $B$ and $\mathbf{x}$ is a vector with $n$ entries from $X$. Then\begin{gather*}\mathcal{H}(n, \mathbf{b}, \mathbf{x})=\frac{1}{q}\binom{t}{n}\left((kq+1)^n+q-1\right).\end{gather*}
\begin{proof}
It is sufficient to count the number of possible choices of $b_i\in B$ and $x_i\in X$ such that \begin{gather*}b_1x_1+\ldots+b_nx_n\equiv0\pmod{q}.\end{gather*} First, we briefly estimate the number of solutions of the congruence $\mathbf{b}\cdot \mathbf{x}\equiv0\pmod{q}$ where $x_i\in X\setminus\{0\}$. We will denote this by $\mathcal{H}^*(n, \mathbf{b}, \mathbf{x})$. For a fixed choice of $\mathbf{b}=\left[b_1, \ldots, b_n\right]$, there are $(kq)^{n-1}$ possible choices for $x_1, \ldots, x_{n-1}$ hence, $(kq)^{n-1}$ sums of the form $b_1x_1+\ldots+b_{n-1}x_{n-1}$. We observe that if \begin{gather*}b_1x_1+\ldots+b_{n-1}x_{n-1}\equiv \beta\pmod{q},\end{gather*} then it suffices to count the number of $x_n\in X$ such that $b_nx_n\equiv -\beta \pmod{q}$. Since $b_n\in \mathbb{Z}_q^{\times}$ this is equivalent to $x_n\equiv-\beta b_n^{-1}\pmod{q}$. We notice that in $X\setminus\{0\}$ every residue class appears exactly $\frac{kq}{q}=k$ times, therefore there are $k$ possible choices for $x_n$. Thus, there are $(kq)^{n-1}\cdot k=\frac{(kq)^n}{q}$ solutions of the congruence $\mathbf{b}\cdot \mathbf{x}\equiv0\pmod{q}$. Since we can choose $\mathbf{b}$ in $\binom{t}{n}$ ways, we conclude that \begin{gather*}\mathcal{H}^*(n, \mathbf{b}, \mathbf{x})=\binom{t}{n}\frac{(kq)^n}{q}.\end{gather*} It suffices to consider the case where $\mathbf{x}$ contains some $0$s. It is evident that if $0$ appears $m<n$ times as an element of $\mathbf{x}$, the number of solutions of $\mathbf{b}\cdot \mathbf{x}\equiv0\pmod{q}$ is equal to \begin{gather*}\mathcal{H}^*(n-m, \mathbf{b}, \mathbf{x})=\binom{t}{n}\frac{(kq)^{n-m}}{q}.\end{gather*} It is straightforward to see that $0$ can appear $m$ times as an element of $\mathbf{x}$ in $\binom{n}{m}$ possible ways. Also, if all the entries of $\mathbf{x}$ are $0$, there are $\binom{t}{n}$ solutions to the congruence $\mathbf{b}\cdot \mathbf{x}\equiv0\pmod{q}$. Thereby, we obtain \begin{gather*}\mathcal{H}(n, \mathbf{b}, \mathbf{x})=\binom{t}{n}+\sum_{m=0}^{n-1}\binom{n}{m}\mathcal{H}^*(n-m, \mathbf{b}, \mathbf{x}).\end{gather*} If we substitute $\mathcal{H}^*(n-m, \mathbf{b}, \mathbf{x})$ from above this is equal to \begin{gather*}\binom{t}{n}\left(1+\sum_{m=0}^{n-1}\binom{n}{m}\frac{(kq)^{n-m}}{q}\right)=\frac{1}{q}\binom{t}{n}\left(q-1+\sum_{m=0}^{n}\binom{n}{m}(kq)^{n-m}\right).\end{gather*} Using the binomial theorem we see that $\sum_{m=0}^{n}\binom{n}{m}(kq)^{n-m}=(kq+1)^n$. Therefore, we proved that \begin{gather*}\mathcal{H}(n, \mathbf{b}, \mathbf{x})=\frac{1}{q}\binom{t}{n}\left((kq+1)^n+q-1\right).\end{gather*}
\end{proof}
\end{lemma}
 With these auxiliary computations in place, we are now ready to prove our last result. 
\begin{theorem}\label{count}
For an odd prime power $q$, let $\mathcal{E}(q)$ be the number of sets  $A\subseteq\mathbb{Z}_q^{\times}$ having non-empty subset sums equidistributed $\pmod{q}$. Then \begin{gather*}\mathcal{E}(q)=\begin{cases} 
      -1+\frac{1}{q}\left((2^r+2)^{\varphi(q)/2r}+(q-1)\cdot 3^{\varphi(q)/2r}\right) & \quad\text{if}\quad  r\equiv 1\pmod{2} \\
       (2^{r/2}+1)^{\varphi(q)/r}-1 &  \quad\text{if}\quad r\equiv 0\pmod{2}.\\
   \end{cases}\end{gather*}
\begin{proof}
It is evident that in Lemma \ref{lemma} we count the number of solutions of congruence \ref{eq:6},with $kq=\frac{2^r-1}{q}$, if all leaders $b_i$ are different. However, it is possible that $A$ is a union of sets with the same leader. 
\begin{case 1*} \emph{The order of $2$ modulo $q$ is odd.}\\
Observe that for every possible choice of ``$+$'' or ``$-$'', if $B_i$ and $B_{i'}$ have the same leader then \begin{gather*}B_i\cup B_{i'}=\{b_i\cdot (-2^{j-1})\}\cup \{b_i\cdot (+2^{j-1})\} :\quad 1\leq j\leq r.\end{gather*} From this, it is evident that the same leader may appear only twice and the sum of the elements of $B_i\cup B_{i'}$ is $0$ modulo $q$. Therefore, it does not affect the condition $\sum_{a\in A}a\equiv0\pmod{q}$.
We estimate $\mathcal{E}(q)$, by focusing on the possible number of pairs of equal leaders. Suppose a union contains $v$ pairs of sets with equal leaders for $v=0,\ldots, t$ and $z$ sets with different leaders for $z=0,\ldots ,t-v$. Notice that the case $v=z=0$ will be excluded at the end, since it represents the empty union, hence the empty set. For $z$ different leaders we may apply Lemma \ref{lemma} with $t=t-v, n=z, kq=2^r-1$ and obtain the value \begin{gather*}\frac{1}{q}\binom{t-v}{z}\left((2^r)^z+q-1\right).\end{gather*}However, there are $\binom{t}{v}$ ways to obtain $v$ pairs of equal leaders. Thus, \begin{gather*}\mathcal{E}(q)=-1+\sum_{v=0}^{t}\binom{t}{v}\sum_{z=0}^{t-v}\frac{1}{q}\binom{t-v}{z}\left((2^r)^z+q-1\right).\end{gather*} From the binomial theorem, we see that 
\begin{gather*}\sum_{z=0}^{t-v}\frac{1}{q}\binom{t-v}{z}\left((2^r)^z+q-1\right)=\frac{1}{q}\left((2^r+1)^{t-v}+(q-1)\cdot 2^{t-v}\right).\end{gather*}This yields \begin{gather*}\mathcal{E}(q)=-1+\sum_{v=0}^{t}\binom{t}{v}\frac{1}{q}\left((2^r+1)^{t-v}+(q-1)\cdot 2^{t-v}\right)\end{gather*} which again, by the binomial theorem is equal to \begin{gather*}-1+\frac{1}{q}\left((2^r+2)^t+(q-1)\cdot 3^t\right).\end{gather*} Since the number of leaders is $t=\frac{\varphi(q)}{2r}$, we obtain \begin{gather*}\mathcal{E}(q)=-1+\frac{1}{q}\left((2^r+2)^{\varphi(q)/2r}+(q-1)\cdot 3^{\varphi(q)/2r}\right).\end{gather*}
\end{case 1*}
 \begin{case 2*} \emph{The order of $2$ modulo $q$ is even.}\\ 
Recall that from the previous analysis if $r$ is even, both $b_i\cdot 2^{j-1}$ and $b_i\cdot2^{j-1+r/2}$ must have the same sign in $B_i$ for every $j=1, \ldots, r/2$, since $2^{r/2}\equiv-1\pmod{q}$. This implies that the sum of the elements of $B_i$ (and thus, of $A$) give a sum equal to $0\pmod{q}$ and that no pair of equal leaders appears. Hence, the number of possible choices for ``$+$'' in $B_i$ is equal to the number of positive integers which have a binary representation consisting of powers of $2$ not greater than $2^{r/2-1}$. This number is $2^{r/2}$. (Again, we want to take account the case where all signs are ``$-$'', which corresponds to $S_i^{+}=0$). For $z$ different leaders, there are $(2^{r/2})^z$ choices of ``$+$'' and $\binom{t}{z}$ ways to choose the leaders. Thus, there are $\binom{t}{z}(2^{r/2})^z$ possible unions. Hence, \begin{gather*}\mathcal{E}(q)=\sum_{z=1}^{t}\binom{t}{z}(2^{r/2})^z.\end{gather*} From the binomial theorem, this is equal to $(2^{r/2}+1)^t-1$. Since in this case, $t=\varphi(q)/r$, we obtain that \begin{gather*}\mathcal{E}(q)= (2^{r/2}+1)^{\varphi(q)/r}-1.\end{gather*}
This completes the proof and the paper.
\end{case 2*}
\end{proof}
\end{theorem}
\makeatletter
\renewcommand{\@biblabel}[1]{[#1]\hfill}
\makeatother

\end{document}